\documentclass[12pt,letterpaper]{article}
\usepackage{t1enc}
\usepackage[latin1]{inputenc}
\usepackage[english]{babel}
\begin{document}
\date{}
\title{Periodic Perturbations  of Non-Conservative Second Order Differential Equations}
\author{A.  Raouf  Chouikha \footnote
{Universite Paris 13 LAGA UMR 7539 Villetaneuse 93430}
}
\maketitle

\begin{abstract}
Let the Lienard system \ $ u'' + f(u) u' + g(u) = 0$ \  with an isolated periodic solution.
This paper concerns the behavior of periodic solutions of Lienard system under small periodic perturbations. \footnote
{2000 Mathematics Subject Classification  \ 34C25, 34C35}

\end{abstract}

 \maketitle
\vspace{1cm}

\section{Introduction}
Consider the second order differential equation of the type $$(E_\epsilon )  \hspace{2cm} x'' + g(t,x,x',\epsilon) = 0$$ where $\epsilon > 0$ is a small parameter, \ $g$ \ is a $T$-periodic function in $t$ and\ $g(t,x,0,0) = g(x)$\ is independent of $t$. \\ The existence problem of non constant periodic solutions of this equation in the case where $g$ is independent on $x'$ and is continuously differentiable has been studied by many authors. \\ Indeed,  in the latter case certain among them  proved existence of solutions of \quad $x'' + g(t,x,\epsilon ) = 0.$\ For a review see Chow-Hale [C-H] and Hale [H].  \\ 
Example given by Hartman proved  non existence cases     
 of that equation if we do not suppose $g$ independent on $x'$.  We then cannot generalize their result.\\ 

Let the following equation, which is a perturbation of Lienard type
$$
(1_\epsilon)
 \hspace{2cm} x'' + f(x) x' + g(x) = \epsilon h(\frac {t}{T},x,x',\epsilon )  
$$
where $h$ is $T$-periodic in $t$, $f$ and $g$ are functions only dependent on $x$, satisfying conditions defined below.  We look for periodic solutions of ($1_\epsilon$) for $\epsilon $ small enough under some additional hypothesis. It is assumed that the unperturbed system has an isolated periodic solution. The perturbation is supposed to be {\it controllably periodic} in the Farkas sense [F], i.e. it is periodic with a period which can be chosen appropriately. \\  We prove an existence theorem for this equation. \\ 

 Loud [L] already proves for $f(x) \equiv c$, the existence of periodic solution of the equation
\begin{equation}
x'' + c  x' + g(x) = \epsilon h(t),
\end{equation}
where the perturbation does not depend on the state.
 He uses for that a variant of the implicit function theorem. More exactly, he considers a function \ $ g(x) = x k(x)$ \ where $k$ is continuously differentiable \ $k(x) > 0, \quad x\neq 0$. \ And either $$ x \frac{d}{dx} k(x) > 0, \quad x\neq 0$$ or $$ x \frac{d}{dx} k(x) < 0 , \quad x \neq 0$$ always holds with the possible exception of isolated points. 
Let us notice  that these conditions imply on one hand the monotonicity of the period function $T$ for the system $x'' + g(x) = 0$.  If $g$ is differentiable they imply on the other hand \ $g''(0) = 0 ,$ \ that is a necessary condition. \\
Let \ $u(t)$ \ be a non-constant $\omega $-periodic solution of the  equation $$ x'' + g(x) = 0$$ and define $$F(s) = \int_0^\infty  u'(t+s) f(t) dt.$$  [L]  observes that if for some $s_0,\quad F(s_0) = 0$ \  while \ $F'(s_0) \neq 0$ \ then for sufficiently small \ $\epsilon > 0$ \ there exists an \ $\omega $-periodic solution \ $v(t,\epsilon )$ \ of the perturbed equation 
 $$ x'' + g(x) = \epsilon f(t) = \epsilon f(t+\omega ) $$

\section{Existence and non-existence of periodic solutions of $E_\epsilon $}

\subsection{A non existence result}
On the other hand, according to P. Hartman ([H],
p. 39), equation ($1_\epsilon $) in general does not have a non constant periodic solution, even if $x g(t,x,x') > 0$. \\
 The following example given by Moser proves the non existence of a non constant periodic solution of 
$$x'' + \phi (t,x,x') = 0.$$  Let  $$\phi (t,x,y) = x + x^3 + \epsilon f(t,x,y), \qquad \epsilon > 0$$  satisfying the following conditions for\quad
$\phi \in C^1(R^3), \quad f(t+1,x,y) = f(t,x,y),$ \quad  with $$f(0,0,0) = 0, \qquad f(t,x,y) = 0  \  {\it if} \  xy = 0$$ $$\frac{\phi }{x} \rightarrow \infty \quad {\it when } \quad x \rightarrow \infty $$ uniformly in $(t,y) \in R^2,$   $$\frac{\delta f   }{\delta y} > 0  \ {\it if}  \ xy >0, \quad {\it and} \quad  \frac{\delta f   }{\delta y} = 0 \ {\it  otherwise}.$$ \ $x,y$\ verifying \ $ \mid x \mid < \epsilon , \mid y \mid < \epsilon .$\\  In fact, we have \ $x f(t,x,y) $\  and $ y f(t,x,y)) > 0 $ \  if $ x y > 0, \mid x\mid < \epsilon , y $\  arbitrary  and $ \phi  = 0 $ otherwise.\\ Notice that \ $\frac{\delta f   }{\delta y} $ \ is small.\\ 

  The function \ $V = 2 x^2 + x^4 + 2 {x'}^2$ \ satisfies \  $V' = -4 \epsilon x' f(t,x,x')$, so that \ $ V' < 0 $ if $ x x' > 0, \mid x \mid < \epsilon $\ and $V' =0$\ otherwise. \\ Thus \ $x$\ cannot be periodic unless \ $V' = 0$   .\\ 

This example is significant because it shows in particular encountered difficulties  in order to establish existence results of periodic solutions for example for the perturbed Lienard equation  \ $x'' + f(x) x' + g(x) = \epsilon h(t,x,x',\epsilon ).$ \ The period of the perturbed equation should be 'controlled' in order to state existence of periodic solution.\\

\subsection{Case where $g$ is independent on $x'$}
Consider equation of the type
\begin{equation}
x'' + \phi (t,x,\epsilon) = 0,
\end{equation}

where $\epsilon > 0$ is a small parameter, $\phi $ is a continuous
function, $T-$periodic on $t$ such that $\phi (t,x,0) = {\tilde g}(x)$.

 More precisely, under the following  hypotheses for the function $g$ defined on $R \times (\alpha ,\beta )\times ]0,\epsilon _0].$ 
\begin{eqnarray}
\cases{
(1) \quad \phi  \ {\it is \ T \ - \ periodic \ on \ t}  & \cr
(2) \quad \phi (t,x,0) = {\tilde g}(x) & \cr
(3) \quad if  \quad  x\neq 0, {\it we \ have} \ {\tilde g}(x) x > 0. & \cr}
\end{eqnarray}

That means for $ \epsilon  = 0$ the autonomous system
\begin{eqnarray}
\cases{
x' = y  & \cr
y' = -  {\tilde g}(x) & \cr}
\end{eqnarray}

has the origin $(0,0)$ as a center.\\ In using a version of fixed
point theorem due to W. Ding, P. Buttazzoni and A. Fonda [B-F]
proved there are periodic solutions of (2) provide that the period \ $T$ \ of the autonomous associated system is monotone and 
$\epsilon $ is small enough. \ $\phi$\ is assumed to be (only) continuous\\ More exactly, they show that if the function $\phi (t,x,\epsilon )$ is continuous, then the
periodic solutions of such equations may be located near to the
solutions of the autonomous equation, provided that periodic
solutions of (4) exist and the period function is strictly
monotone. Moreover, there is a solution making exactly $N$ rotations
around the origin in the time $kT$. \\
 Their results improve those of Loud [L] , that thought that the function $\phi $ had to be continuously differentiable.\\
Moreover, under the light of the preceding example [H], it seems that the methods described above does not generalize if one supposes  $ \phi$ dependent on $x'$  : \ $\phi \equiv \phi (t,x,x',\epsilon )$.\\ So, an other condition on the period appears to be necessary to obtain existence of periodic solutions of the perturbed equation.\\
One nevertheless can show an analogous result to the preceding one under more restrictive hypotheses.  In particular on the appropriated choice of the period of a periodic solution of the perturbed Lienard equation.\\

\section{A controllably periodic perturbation}

One refers to a method due to Farkas inspired of the one of Poincaré.  The determination of controllably periodic perturbed solution. This method proved to be itself very effective particularly for the perturbations of various autonomous systems. Since we know it for example a good application for perturbed Van der Pol  equations type [F-F].  We logically may expect that the Farkas method be again applied for perturbed  Lienard equations.
This has been considered and proved by Farkas himself  [F] .\\ Our proof we give here made it more simple and contains some modifications in using in particular methods of [F1] to estimate existence regions of periodic solutions for that equation. \\
The perturbation is supposed to be `controllably periodic', i.e., it is periodic with a period which can
be chosen appropriately. Under very mild conditions it is proved that to each small enough amplitude of the perturbation
there belongs a one parameter family of periods such that the perturbed system has a unique periodic solution with this
period.\\

\subsection{Basic hypotheses}
Let us consider the (unperturbed) Lienard equation $$(L)\qquad u'' + f(u) u' + g(u) = 0.$$
In order to have an unique periodic solution we have to suppose the functions \ $f$\ and\ $g$\ are of class \ $C^2$.\ $f$ is even, and \ $g$\ is odd. The integral $$F(x) \int_0^x h(t) dt ,\qquad G(x) \int_0^x h(t) dt$$
  of \ $f$\ and\ $g$\ respectively are such that \ $ lim_{x\rightarrow \infty} F(x) = \infty,$\ and \ $ lim_{x\rightarrow \infty} G(x) = \infty.$\ It is assumed that \ $F$ \ has a unique zero.
Then it is known [C-L] that (L) has a stable non constant periodic solution \ $u_0(t)$\ with period \ $\tau _0$.\\ 
Equation (L) is usually studied by means of an equivalent plane system. The most used ones is :

\begin{eqnarray} 
\cases{
u' = v  & \cr 
v' = - g(u) - f(u) v & \cr}  
\end{eqnarray}
 
and also

\begin{eqnarray} 
\cases{
u' = v - F(u) & \cr 
v' = - g(u) & \cr}  
\end{eqnarray}

In fact, there are equivalent to the 2-dimensional system 
$$(S) \qquad \dot x = h(x)$$
after introducing the notations
 \ $x = col [x_1,x_2]$

\begin{eqnarray} 
\cases{
x_1 = - \dot u(t) - F(u(t)) & \cr 
x_2 = u(t) & \cr}  
\end{eqnarray}
 
where \ $x = col [x_1,x_2]$\ and \ $h(x) = col [g(x_2) , -x_1-F(x_2(t)).$\\ 
Suppose $$ u_0(0) = a, \qquad u'(0) = 0 > 0$$ so that   
the periodic solution of period $\tau _0$ of the variational system
$$\dot y = {h'}_x(p(t)) y$$ 
 is 
$$p(t) = col [-\dot u_0(t)-F(u_0(t)),  u_0(t) ].$$ 
We then have
$$ \dot p(t) = col[g(u_0(t)), \dot u_0(t)].$$ 
So,  the initial conditions are $$ p(0) = col[-F(a),a], \qquad  \dot p(0) = col[g(a), 0].$$

\bigskip

Let us consider the following perturbed Lienard equation of the form
$$(L_R)\qquad \ddot u + f(u) \dot u + g(u) = \epsilon \gamma (\frac {t}{\tau }, u, \dot u)$$ where \ $t \in \!R,\quad \epsilon \in \!R$\ is a small parameter, \ $\mid \epsilon \mid < \epsilon _0,\quad \tau $\ is a real parameter such that \ $\mid \tau  - \tau _0\mid < \tau _1$\ for some \ $0 < \tau _1 < \frac {\tau _0}{2}$.\\
Moreover, the closed orbit $$\{ (u,v) \in \!R^2 : u(t) = u_0(t), \ v(t) = \dot u_0(t), \ t \in [0,\tau _0]\}$$
belongs the region \ $\{ (u,v) \in \!R^2 : u^2 + v^2 < r^2 \}.$  \ $q $\ is a function of class \ $C^2,$\ $\tau $-periodic in $t$.\\
By the same way as for $(L)$, the 2-dimensional equivalent system to $(L_R)$ is 
$$(S_L)\qquad \dot x = h(x) + \epsilon q(\frac {t}{\tau }, x)$$
where \ $q = col [ q_1,q_2],$
  
\begin{eqnarray} 
\cases{
q_1 = - \gamma (\frac {t}{\tau }, x_2, -x_1 - F(x_2))  & \cr 
q_2 = 0 & \cr}  
\end{eqnarray}

\subsection{Existence of periodic solutions of $(L_R)$}

Now we will use Poincare method for the determination of the approximated solution of the perturbed equation $(L_R)$.\ It is assumed the existence of the fondamental matrix solution of the first variational system of  \ $\dot x = h(x)$ \ and the unique periodic solution \ $p(t)$\ corresponding to \ $u_0(t).$\\  
In order to get estimates for the existence of periodic solutions we have to calculate some constants.
Following Farkas [F], the Jacobi matrix \ $J$\ has the following form 
$$J(\tau _0) = - I + \pmatrix{g(a) & 0\cr 0 & 0 \cr} + Y(\tau _0) $$ 
$I = Id_2$ \ and \ $Y(t)$\  is the fundamental solution matrix of the varational system with \ $Y(0) = I$

\begin{equation} 
 \dot y = \pmatrix{O & g'(u_0(t)) \cr - 1 & - f(u_0(t))} y.
\end{equation}

  It is proved, that if \ $det J(\tau _0) \neq 0$\  then there exist uniquely determined functions \ $\tau (\epsilon ,\phi )$\ and \ $h(\epsilon , \phi )$\ defined in the neighborhood of \ $(0,0)$\ such that the function $$u(t;\phi, p_0 + h(\epsilon , \phi ), \epsilon , \tau (\epsilon ,\phi ))$$ is periodic solution of system $(S_L)$ so that \ $\tau (0,0) = \tau _0, \ h(0,0) = 0$.\  Moreover, it is given an estimate for the region in which the the variables \ $\epsilon $\ and \ $\phi $\ may vary . It is needed for that to evaluate the norm of the difference of Jacobi matrices\ $J(\epsilon ,\phi ,\tau ,h) -    J(0,0,\tau _0,0).$\\

A trivial $\tau _0$-periodic solution of (9) is \ $col [g(u_0(t)), \dot u_0(t)].$\ A trivial calculation gives the other linearly independent solution of (9)
$$col [g(u_0(t)) v(t), \dot u0(t) v(t) + \frac{g(u_0(t))}{g'(u_0(t))} \dot v(t)] $$
where 
$$v(t) = \int_0^t [g_0(s)]^{-2} g'(u_0(t)) exp [-\int_0^s f(u_0(\sigma ))d\sigma ]ds$$
for \ $t \in [0, \tau _0]$.\\

Then the fundamental solution matrix of (9) with \ $Y(0) = I$\ is 
$$Y(t) = \pmatrix{\frac{g(u_0(t))}{g(a)} & g(a)g(u_0(t)) v(t) \cr \frac{\dot u_0(t)}{g(a)} & g(a)\dot u_0(t) v(t) + g(a)\frac{g(u_0(t))}{g'(u_0(t))} \dot v(t)}.$$
According to Liouville's formula  the Wronskian determinant \ $W(t)$\ with \ $W(o) = 1$ is given by 
$$W(t) = exp [- \int_0^t f(u_0(\tau )) d\tau ] .$$ 
The characteristic multipliers of (9) are \ $\rho _1 = 1$\ and   $$\rho _2 = W(\tau _0) = exp [- \int_0^{\tau _0} f(u_0(\tau )) d\tau ] .$$
$\rho _2 < 1$ \ if and only if 

\begin{equation} 
 \int_0^{\tau _0} f(u_0(\tau )) d\tau > 0.
\end{equation}
             
The initial conditions give
$$Y(\tau _0) = \pmatrix{1 & g^2(a) v(\tau _0) \cr 0 & \rho _2}.$$
Thus, we get 
$$J = \pmatrix{g(a) & g^2(a) v(\tau _0) \cr 0 & \rho _2 - 1},$$ 
$$J^{-1} = \pmatrix{g^{-1}(a) & g^2(a) v(\tau _0) (1 - \rho _2)^{-1} \cr 0 & -(1 - \rho _2)^{-1}}.$$
Therefore, $$\mid\mid J^{-1} \mid\mid = 2 \ max \ [g^{-1}(a) , (1 - \rho _2)^{-1} , g^2(a) v(\tau _0) (1 - \rho _2)^{-1} ]. $$ 

The inverse matrix of $Y(t)$ is  
$$Y^{-1}(t) = W(t) \pmatrix{g(a)\dot u_0(t) v(t) + g(a)\frac{g(u_0(t))}{g'(u_0(t))} \dot v(t) & - g(a)g(u_0(t)) v(t) \cr -\frac{\dot u_0(t)}{g(a)} & \frac{g(u_0(t))}{g(a)}}  .$$
  
Now we have to determine the constants for system $(L_R)$. Following [F] let us denote $$S = \{ x=(x_1,x_2) \in \!R^2\  / \  {x_2}^2 + [- x_1 - F(x_2(t))]^2 < r^2 \}$$ 
\begin{eqnarray} 
\cases{
 g_0 := max_{x\in S} \mid g(x_2) \mid , & \cr 
g_1 :=max_{x\in S} \mid g'(x_2) \mid , & \cr
g_2 = max_{x\in S} \mid g''(x_2) \mid . & \cr}  
\end{eqnarray} 

$$f_1 := max_{x\in S} \mid f(x_2) \mid , \qquad f_2 :=max_{x\in S} \mid f'(x_2) \mid $$  
\begin{eqnarray} 
\cases{
 q_0:=max_{x\in S, s\in \!R} \mid q(s,x) \mid, & \cr 
q_1:=max_{x\in S, s\in \!R} \mid q'_x(s,x) \mid, & \cr
q_2:= max_{x\in S, s\in \!R} \mid q'_s(s,x) \mid . & \cr}  
\end{eqnarray} 
$$ K:= max_{t\in [\frac{-\tau _0}{2},\tau _0]} \mid Y(t) \mid, \qquad K_{-1}:= max_{t\in [\frac{-\tau _0}{2},\tau _0]} \mid Y^{-1}(t) \mid .$$
Thus, we may deduce that
$$ P:= max_{t\in [\frac{-\tau _0}{2},\tau _0]} \mid \dot p(t) \mid \leq \frac{K}{2},$$
The phase initial $\phi$ and the period $\tau$ have to verify the following, which can be easely obtained from the above estimates, see [F]
$$\phi < \frac{\tau _0 }{2}, \qquad \mid \tau  - \tau _0 \mid < \frac{\tau _0 }{2}.$$
If in addition we suppose $\epsilon $ and $h$ are such that 
$$\frac {3}{2} g_0 \mid \epsilon \mid + \mid h \mid < \sigma exp(- \frac{3}{2} g_1 \tau _0)$$ 
(here $\sigma $ is the distance between the path of the periodic solution and the boundary of $S$)\ then a solution 
of $(L_R)$ exists.\\
We may resume in the following

\bigskip

{\bf Theorem 1}\qquad {\it If $1$ is a simple characteristic multiplier of (9) (that means inequality (10) holds) then there are two functions \ $\tau , h : U \rightarrow \!R$ \ and a constant \ $\tau _1 < \frac {\tau _0}{2}$ \ such that the solution \ $ u(t,\phi, a+h(\epsilon ,\phi ),\epsilon, \tau ) $\ of equation 
$$(L_R)\qquad \ddot u + f(u) \dot u + g(u) = \epsilon \gamma (\frac {t}{\tau }, u, \dot u)$$
exists for \ $(\epsilon ,\phi )\in U$, and satisfying properties \ $ \mid \tau - \tau _0 \mid < \tau _1,\ \tau (0,0) = \tau _0, \ h(0,0) = 0 .$}

\bigskip

\subsection{Special cases}
As a corrolary we may deduce from the above some results about autonomous perturbations of the Lienard system
$$(L_{R_A}) \qquad \ddot u + f(u) \dot u + g(u) = \epsilon \gamma $$
where the perturbation is independent on the time variable \ $\gamma \equiv \gamma (u, \dot u, \epsilon , \tau )$.\\ The equivalent plane system is of the form
$$\dot x = h(x) + \epsilon q(\frac {t}{\tau }, x)$$
where \ $q = col [ q_1,q_2].$\\
Consider the following autonomous perturbed Lienard equation of the form
$$(L_{R_A})\qquad \ddot u + f(u) \dot u + g(u) = \epsilon \gamma \equiv \gamma (u, \dot u, \epsilon , \tau )$$ where \ $t \in \!R,\quad \epsilon \in \!R$\ is a small parameter, \ $\mid \epsilon \mid < \epsilon _0,\quad \tau $\ is a real parameter such that \ $\mid \tau  - \tau _0\mid < \tau _1$\ for some \ $0 < \tau _1 < \frac {\tau _0}{2}$.\\
Moreover, the closed orbit $$\{ (u,v) \in \!R^2 : u(t) = u_0(t), \ v(t) = \dot u_0(t), \ t \in [0,\tau _0]\}$$
belongs the region \ $S =\{ (u,v) \in \!R^2 : u^2 + v^2 < r^2 \}.$  \ $\gamma $\ is a function of class \ $C^2.$ \\ In this case the perturbation is independent on the initial phase \ $\phi $.\\
We then have the following

\bigskip

{\bf Corrolary 2}\qquad {\it Suppose inequality (10) to hold, then there are constants \ $\epsilon  _0$ \ and  \ $\tau _1 < \frac {\tau _0}{2}$ \
such that to each \ $\epsilon \in [-\epsilon _0,\epsilon _0]$\ there exist two functions \ $\tau , h $ \ only dependent on \ $\epsilon : \tau \equiv \tau(\epsilon ),\ h \equiv h(\epsilon )$ \ such that the equation 
$$(L_{R_A})\qquad \ddot u + f(u) \dot u + g(u) = \epsilon \gamma \equiv \epsilon \gamma (u, \dot u, \epsilon , \tau )$$
 has a unique periodic non constant 
solution \ $u(t, \epsilon )$\ with period \ $\tau $.\ Moreover,\ $\tau (0) = \tau _0$ , and satisfy the properties \ $ \mid \tau - \tau _0 \mid < \tau _1,\ \tau (0) = \tau _0, \ h(0) = 0 .$}

\bigskip

On other hand,  a second special case may occur when the perturbation does not depend on the state of the system. That means the perturbation is independent on $u$ 
$$\gamma \equiv \gamma (\frac {t}{\tau }).$$
Then the above estimates can be easely calculated.

\vspace{1cm}

\begin{center}
\Large
REFERENCES\\
\end{center}

\bigskip

[B-F] \qquad P.Buttazzoni and A.Fonda \quad {\it Periodic perturbations of scalar second order differential equations} \quad Discr. and Cont. Dyn. Syst., vol 3, n° 3, p. 451-455, (1997).\\

[C-H]\qquad S.N.Chow-J.Hale \quad {\it Methods of bifurcation theory}, \quad Springer, Berlin (1982).\\

[C-L]\qquad E.Coddington-N.Levinson \quad {\it Theory of ordinary differential equations},
\quad Mc Graw-Hill, New-York, (1955).\\

[F-F] \qquad I.Farkas-M.Farkas \quad {\it On perturbations of Van der Pol's equation} \quad Ann. Univ. Sci. Budapest, Sect math., p. 155-164, (1972).\\

[F] \qquad M.Farkas \quad {\it Periodic perturbations of autonomous systems}\quad Alkalmaz. Math. Lapok, vol 1, n° 3-4, p. 197-254, (1975).\\

[F1] \qquad M.Farkas \quad {\it Estimates on the existence regions of perturbed periodic solutions},\quad SIAM J. Math. Anal., vol 9, p. 867-890, (1978).\\

[Ha]\qquad J.Hale \quad {\it Topics in dynamic bifurcation theory}, \quad Wiley, New-York, (1981).\\

[H] \qquad P.Hartman \quad {\it On boundary value problems for superlinear second order differential equation}\quad J. of Diff Eq. , vol 26, p. 37-53, (1977).\\

[Lo] \qquad W.S. Loud \quad {\it  Periodic solutions of $x'' + c x' + g(x) = \epsilon f(t)$  }
\quad Mem. Amer. Math. Soc., n 31, p. 1-57, (1959).\\

\end{document}